\documentclass[12pt]{article}
\usepackage{bbm,pifont,eufrak,eucal,mathrsfs,latexsym,graphicx,epsfig,color,rotating,multirow}

\usepackage[ansinew]{inputenc}

\newtheorem{thm}{Theorem}[section]

\newtheorem{cor}[thm]{Corollary}

\newtheorem{defi}{Definition}

\newtheorem{lem}[thm]{Lemma}

\newtheorem{pro}[thm]{Proposition}

\newcommand{\ad}{\mathrm{ad}}
\newcommand{\Aut}{\mathrm{Aut}}
\newcommand{\cf}{\emph{cf }}
\newcommand{\End}{\mathrm{End}}
\newcommand{\equi}{\Longleftrightarrow}
\newcommand{\Frac}{\mathrm{Frac}}
\newcommand{\hcf}{\mathrm{hcf}}
\newcommand{\Id}{\mathrm{Id}}
\newcommand{\ie}{\emph{i.e. }}
\newcommand{\inv}{^{-1}}
\newcommand{\kk}{{\mathbbm k }}
\newcommand{\mod}{\mathrm{ \: mod \:}}
\newcommand{\qq}{\forall\;}
\newcommand{\resp}{\emph{resp. }}
\newcommand{\sub}{\subseteq}

\newcommand{\tq}{\, : \,}
\newenvironment{dem}[1][]{%
{\it Proof #1 : }\it}{%
\hspace*{\fill}\nolinebreak[1]\hspace*{\fill} {\it q.e.d.}\\}

\begin{document}

\title{On the automorphism group of the first Weyl algebra}
\author{M.K. Kouakou and A. Tchoudjem \\
Universit\'{e} de Cocody\\
UFR-Math\'{e}matiques et Informatique\\
22 BP 582 Abidjan 22 C\^{o}te d'Ivoire\\
e-mail : makonankouakou@yahoo.fr\\
and\\
Université Claude Bernard Lyon I\\
Institut Camille-Jordan\\
69622 Villeurbanne cedex France\\ 
tchoudjem@math.univ-lyon1.fr}
\date{\today}
\maketitle

\begin{abstract}

Let $A_{1} := \kk [t, \partial ]$ be the first algebra over a field $\kk$ of characteristic
zero. One can associate to each right ideal $I$ of $A_1$ its Stafford subgroup, which is a subgroup of $\Aut_\kk(A_1)$,  the automorphism group of the ring $A_1$. In this article we show that each Stafford subgroup is equal to its normalizer. For that, we study when the Stafford subgroup of a right ideal of $A_1$ contains a given Stafford subgroup.

\end{abstract}

\section*{Introduction}
Let $\kk$ be a commutative field of characteristic zero. We note $A_{1}$ the first Weyl agebra over $\kk$ \ie :\[A_1:=A_{1}(\kk)=\kk[t,\partial]\]
where $\partial,t$ are related by $\partial t -t\partial=1$.

\begin{defi}

For  a right ideal $I$ of $A_{1}$, the Stafford subgroup associated to $I$  is :
$$H(I) := \left\{\sigma \in Aut_{\kk}(A_{1}) : \sigma(I) \simeq I\right\}$$
(where the symbol ``$\simeq$'' means ``$\sigma(I)$ is isomorphic to $I$ as a right-$A_1-$module'').

\end{defi}

By \cite{sta}, it is known that each subgroup $H(I)$ is isomorphic
to an automorphism group $\Aut_{\kk}(\mathcal{D}(X))$, where $\mathcal{D}(X)$ is the $\kk-$algebra
of differential operators over an algebraic affine curve $X$.

A natural question is : 

``are the Stafford subgroups normal in $\Aut_\kk(A_1)$'' ? 

The answer is no.

Stafford showed that if $X_2$ is the famous algebraic plane curve defined by the equation :
\[x^2=y^3\]
and if $I_2$ is the right ideal of $A_1$ :
\[I_2:=\left\{d \in A_1 \tq  d(\kk[t])  \sub \kk[t^2,t^3] \right\}\]
then the subgroup $H(I_2)$ is isomorphic to $\Aut_{\kk}(\mathcal{D}(X_{2}))$ and is equal to its own normalizer
in $ \Aut_{\kk}(A_{1})$.

We will show in this paper that the subgroup $H(I)$ is equal to
its own normalizer for all right ideal $I$ of $A_1$.

We begin by  giving some definitions and by fixing  some notations that will be used in
this paper.



\section{Definitions and some properties }

The ring $A_{1}$ contains the subrings $R:=\kk[t]$ and $S:=\kk[\partial]$. It is well known that $A_{1}$ is a two-sided noetherian integral domain. Since the characteristic of $\kk$ is zero, $A_{1}$ is also hereditary (\cf \cite{ss}) \ie every non zero right ideal of $A_1$ is a projective right-$A_1-$module. 

The ring $A_1$ has a quotient divison ring, denoted by $Q_{1}$. For any finitely generated right-submodule $M$ of $Q_{1}$, the dual $M^{*}$, as a left-$A_{1}-$module will be identified with the set $\{u\in Q_{1}:uM\subseteq A_{1}\}$, and $End_{A_{1}}(M)$ with
 the set $\{d\in Q_{1}:dM\subseteq M\}$ (\cf \cite{sta}).

The division ring $Q_{1}$ contains the subrings $D:=\kk(t)[\partial]$ and $E:=\kk(\partial)[t]$. The elements of $D$ are $\kk$-linear endomorphisms of $\kk(t)$. More precisely, if

$d:=a_{n}\partial^{n}+...+a_{1}\partial+a_{0}$ for some  $a_{i}\in \kk(t)$ and if $h\in \kk(t)$, then :
\[
d(h):= a_{n}h^{(n)}+... +a_{1}h^{(1)}+a_{0}h \;,
\]

where $h^{(i)}$ denotes the $i$-th derivative of $h$ and $a_{i}h^{(i)}$ is a product in $\kk(t)$. We note that:
\[\qq d,\, d' \in \kk(t)[\partial],\,\qq h \in \kk(t),\,   (dd')(h)=d(d'(h))  \;.\]

For $V$ and $W$ two vector subspaces of $\kk(t)$, we set :
\[\mathcal{D}(V,W):=\{d\in \kk(t)[\partial]:d(V)\subseteq W\}\;.\]

Notice that $\mathcal{D}(R,V)$ is a right  $A_{1}-$submodule of $Q_{1}$ and $\mathcal{D}(V,R)$ is a left $A_{1}-$submodule of $Q_{1}$. If moreover $V \subseteq R$, then  $\mathcal{D}(R,V)$ is a right ideal of $A_{1}$. When $V=R$, one has $\mathcal{D}(R,R)=A_{1}$.

If $I$ is a right ideal of $A_{1}$, we set  : \[I \star1:=\{d(1) \tq d\in I\}\;.\]

It is clear $I\star1$ is a $\kk-$vector subspace of $\kk[t]$ and that : \[I\subseteq \mathcal{D}(R,I\star1)\;.\]

The inclusions $A_{1}\subset k(\partial)[t]$ and $A_{1}\subset k(t)[\partial]$ show that, at least, two notions of degree can be defined on $A_{1}$ : the degree in ``$t$" or $t-$ degree and the degree in ``$\partial$" or $\partial-$degree. Naturally, those degree notions extend to $Q_{1}$. We will note them, respectively, $\deg_t$ and $\deg_\partial$.

\bigskip

\section{Primary decomposable subspaces}

In order to describe the right ideals of $A_1$, it is convenient to use the notion of  
{\it primary decomposable subspaces} of $\kk[t]$.

Recall that $\kk$ is not necessarily algebraically closed.

Let $b, h \in R = \kk[t]$ and $V$ a $k$-subspace of $R$. We set:
\[
 {\cal O}(b) := \{a \in R : a'\in bR\} \,,\]
where $a'$ denotes the formal derivative of $a$. 

{\it E.g. :} one has $\mathcal{O}(t^{n-1})=\kk + t^n\kk[t]$.

We set also :
\[
S(V ) := \{a \in R : aV \subseteq V \}\mbox{ and } {\cal C}(R, V ) := \{a \in R : aR \subseteq V \}\;.\]

Clearly ${\cal O}(b)$ and $S(V )$ are $\kk$-subalgebras of $R$. If $b\neq 0$, the Krull dimension of $\mathcal{O}(b)$ is $\dim_{\kk}({\cal O}(b))=1$. 

The set $\mathcal{C}(R, V )$ is an ideal of $R$
contained in both $S(V )$ and $V$. Moreover, if ${\cal O}(b) \sub S(V)$ then $b^2R \sub V \ie b^2 \in {\cal C}(R,V)$.  

\begin{defi}
 A non-zero $k$-vector subspace $V$ of $k[t]$ is said to be primary decomposable (p.d. for short) if $S(V )$ contains a $k$-subalgebra $\mathcal{O}(b)$, with $b\neq 0$. In this case $\mathcal{C}(R, V )$ is a non zero ideal of $R$. A p.d. subspace $V$ of $k[t]$ is said irreducible (p.d.i.) if $V$ is not contained in a proper ideal of $k[t]$.
\end{defi}

In \cite{caho}, R.C. Cannings and M.P. Holland have shown that for p.d. $V$ of $R$, there is the equality \[\mathcal{D}(V,V) = \End_{A_{1}}(\mathcal{D}(R, V )) \;.\]

It is shown in \cite{kthese} that for any non zero right ideal $I$ of $A_{1}$, there
exists $x \in Q_{1}$ and  $\sigma\in \Aut_{k}(A_{1})$ such that :
\[x\sigma(I) = \mathcal{D}(R,\kk[X_n]))\,,\]
where $n \in \mathbbm{N}$ and $\kk[X_n]: = k + t^{n}k[t]$ is the ring of regular functions on an affine algebraic affine curve $X_{n}$.

We will show that the inclusion :
\[H(\mathcal{D}(R,\kk[X_n])) \subseteq H(\mathcal{D}(R, V ))\] where $V$  is a proper p.d.i. subspace of $R$, implies :
\[\kk[X_n] = V \;.\]

That result will  lead us  to the conclusion that the subgroup $H(\mathcal{D}(R,\kk[X_n]))$ is equal to its own normalizer in $\Aut_{k}(A_{1})$.

\section{The characteristic elements of a right ideal}

The first step in the classification of right ideals of the first Weyl algebra $A_{1}$ is the following :

\begin{thm}[Stafford  {\cite[lemma 4.2]{sta}}]
If  $I$ is  a non-zero right ideal of $A_{1}$, then there exist $e,e' \in Q_{1}$ such that:
\[(i)\; \; eI\subseteq A_{1} \mbox{ and } eI\cap k[t] \neq \{0\}\;;\]
\[ (ii)\; \; e'I\subseteq A_{1} \mbox{ and } e'I\cap k[\partial] \neq \{0\}\;.\] 
\end{thm}

 With $(i)$ we see that any non-zero right ideal $I$ of $A_{1}$ is isomorphic to another ideal $I'$ such that $I'\cap k[t]\neq \{0\}$. 
 \vskip .5cm

{\it Remark:} the element $e$ (\resp $e'$) of the theorem is a minimal $\partial-$degree element of $I^*$ (\resp a minimal $t-$degree element of $I^*$).

\begin{cor}
There exists an unique element (modulo the multiplicative group $\kk^*$) $f\in I$  such that the full set of elements of $I$ with minimum $t-$degree be exactly :
\[f \kk[\partial]\;.\]

In the same way, there exists an unique element (modulo the multiplicative group $\kk^*$) $e^* \in I^*$   such that the full set of elements of $I^*$  with minimum $t-$degree be exactly :
\[\kk[\partial]e^* \;.\]
\end{cor} 
\begin{dem}
For example,  for $f$ : let $e' \in I^*$ such that $e'I \cap \kk[\partial] \neq \{0\}$. Let $s\in \kk[\partial]$ such that $e'I \cap \kk[\partial] =s \kk[\partial]$. We can take : $f:= {e'}\inv s$.
\end{dem}

\begin{defi}
 The elements $e^* \in I^*$ and $f\in I$ are called the characteristic elements of the ideal $I$. 

If $V$ is a $p.d.$ subspace and if $I=\mathcal{D}(R,V)$, we will also say that $e^*$ and $f$ are the characteristic elements of the p.d.i. subspace $V$.

\end{defi}

{\it Remark :} if $e^*,f$ are the characteristic elements of a right ideal $I$, then $e^*f \in \kk[\partial]$ and :
\[e^*I \cap \kk[\partial] = e^*f\kk[\partial] \;.\]

{\it E.g. :} the characteristic elements of $\kk[X_n]=\kk +t^n\kk[t]$ are :
\[e_n^*:= t^{-n}(t\partial) \in \mathcal{D}(\kk[X_n],R) \;,\]
\[ f_n := (t\partial -1) ... (t\partial -(n-1)) \in \mathcal{D}(R,\kk[X_n]\]
and $e_n^*f_n =\partial^n$.

We now recall  some  important properties of  the p.d.i. subspaces $V$ and of  the associated right
ideals $\mathcal{D}(R, V )$. 

\begin{lem}
Let $I$ be a right ideal of $A_{1}$ such that $I \cap \kk[t] \neq \{0\}$.

--- If $V := I \star 1 := \{d(1) : d \in I\}$, then $V$ is a p.d. subspace of $R$ and
\[I = \mathcal{D}(R, V)\;.\]

---  For any p.d. subspace $W$ of $R$, one has
\[
\mathcal{D}(R,W)\star 1 = W \mbox{ and } \mathcal{C}(R,W) = \mathcal{D}(R,W) \cap \kk[t] \;.\]
\end{lem}

For a proof \cf \cite[theorem \S 3.2]{caho}.
\bigskip

Henceforth, $\theta$ will denote the $\kk$-automorphism of $A_{1}$ such that :
\[\theta(\partial) = t \mbox{ and } \theta(t) = -\partial\;.\]

According to \cite{caho}, the above lemma has the following consequence :

\begin{cor}
For any non zero right ideal $I$ of $A_{1}$,

i) there is a unique $x \in Q_{1}$ (modulo the multiplicative group $\kk^*$) and a unique p.d.i. subspace $V$ of $R$ such that $xI = \mathcal{D}(R, V )$ ;

ii) there is a unique $y\in Q_{1}$ (modulo the multiplicative group $\kk^*$)
and a unique p.d.i. subspace $W$ of $R$ such that : $\theta(yI) = \mathcal{D}(R,W)$.
\end{cor}
\bigskip

Now let us give some properties which characterize a p.d.i. subspace $V$ of $R$.

\begin{pro}\label{pro:cinq}
 Let $V$ be a p.d.i. subspace of $R$ and $m := dim_{\kk}{R}/{V}$.

i) For any $0\neq d\in \mathcal{D}(R, V )$, $d(R)\subseteq V$ and  $\deg_{t}(d) \geq dim_{\kk} R/V$.

ii) If $0 \neq f\in \mathcal{D}(R, V )$ has minimal $t$-degree, then $f(R) = V$.

iii) Let $e^*$ and $f$ be the characteristic elements of $V$. As $e^*\in {\cal D}(R,V)^*$ and as ${\cal D}(R,V)\cap \kk[t] \neq \{0\}$, we have $e^* \in \kk(t)[\partial]$. More over if :
\[
f = t^{m}c_{m }(\partial) + t^{m-1}c_{m-1} (\partial) +...+ c_{0 }(\partial)\]
for some  $c_{i} (\partial) \in \kk[\partial]$, and if 
\[
e^* = b_{m}(\partial)t^{-m} + u\]
 where $b_{m}(\partial) \in \kk[\partial]$, $u \in \kk(t)[\partial]$ and $ deg_{t}(u) < -m$,
then :
\[e^*f = b_{m}(\partial)c_{m}(\partial)\;.\]
\end{pro}

Those properties have all been  proved in \cite[remarques 1,2,3]{kou}.

{\it Remarks :}

---  Note  that the $\kk$-vector space $R/V$ has  finite dimension since  $\{0\}\neq \mathcal{C}(R,V)\subseteq V$.

--- How to calculate $f$? We take any $f' \in \mathcal{D} (R, V )$ with minimal $t-$degree $m$, and we  expand $f'$ as polynomial in $t$:
\[f' =t^{m}a_{m }(\partial) + t^{m-1}a_{m-1} (\partial) +... + a_{0 }(\partial)\]
 where $a_{i} (\partial) \in \kk[\partial]$ for all $i$.

If $p(\partial) := \hcf (a_{m}(\partial), a_{m-1}(\partial),..., a_{0}(\partial))$ then we get $f' = f p(\partial)$ (modulo $\kk^*$).

\section{About the automorphisms that stabilizes  an ideal}

Let $I$ be a right ideal of $A_{1}$ such that $I \cap \kk[t] \neq \{0\}$ or $I \cap \kk[\partial] \neq \{0\}$. In
this paragraph we wish to determine the automorphisms $ \sigma\in Aut_{\kk}(A_{1})$ such
that $\sigma(I) = I$. 

We introduce some particular automorphisms of $A_1$.

If $p \in R$ we define $\sigma := \exp(\ad(p))$ by :
\[\qq d \in A_1,\,\sigma(d) := d+[d, p]+ \frac{1}{2!} [[d, p],p]+ \frac{1}{3!} [[[d, p], p], p]+ ...\]
where $[d, p] :=dp-pd$ for all $d \in A_1$.

As the application :
\[A_1 \to A_1\,, \;\; d \mapsto [d,p]\]
is a locally nilpotent derivation, $\sigma$ is a well defined automorphism of the ring $A_1$. Moreover $\sigma\inv =\exp(\ad(-p))$.
 
The following theorem is fundamental in this paper.

\begin{thm}\label{thm:exp}
Let $V$ be a p.d. subspace of $R$ and $\sigma := exp(ad(p))$ where $p \in k[t]$.
Then $\sigma(\mathcal{D}(R, V )) = \mathcal{D}(R, V )$ if and only if $p \in S(V )$.
\end{thm}

\begin{dem}

{\textbf{Suppose that  \mathversion{bold} $ {p \in S(V )}$}}.

 Let $d\in \mathcal{D}(R,V)$.

Clearly $dp$ and $pd$ are both in $\mathcal{D}(R, V )$, so $[d, p]\in \mathcal{D}(R, V )$, and this implies the first inclusion : $\sigma(D(R, V ))\subseteq\mathcal{D}(R, V )$. In the same way, since $-p \in S(V )$ we get the second inclusion $\sigma^{-1}(D(R, V ))\subseteq \mathcal{D}(R, V )$ and then the
equality $\sigma(\mathcal{D}(R, V )) = \mathcal{D}(R, V )$.
\bigskip

{\textbf{Now suppose that \mathversion{bold} $\sigma(\mathcal{D}(R, V ))\subseteq \mathcal{D}(R, V )$.}}

Let us take an element $f\in \mathcal{D}(R, V )$.

In the formal power series ring $\kk[[T]]$, let \[\log(1+T):=\sum_{k=1}^\infty \frac{(-1)^{k-1}}{k}T^k  \;.\]

In the ring $\kk[[T]]$, there is the equality :
\[\log(1+(e^T-1)) = T \;.\]

If we specialize in $T=\ad p$, we get :\[\sum_{k=1}^\infty \frac{(-1)^{k-1}}{k}(\sigma-\Id_{A_1})^k(f) = [f,p]\]
(the left hand side sum is finite because $f \in A_1$ and $p \in R$).

Now, all the terms   $f,\sigma(f), \sigma^{2}(f), ..., \sigma^n(f),...$ belong to $\mathcal{D}(R, V )$, so do the terms :
\[(\sigma-\Id_{A_1})^k(f) \]
($k \ge 1$). Therefore, $[f,p] \in {\cal D}(R,V)$ and we have :
\[[f,p](1) =f(p) -pf(1) \in V\]
\[\Longrightarrow pf(1) \in V \;.\]

But we have :
\[V = {\cal D}(R,V)\star 1\]
\[= \left\{f(1) \tq f \in {\cal D}(R,V)\right\}\]
so $pV \subseteq V$ \ie $p\in S(V)$.

\end{dem}

\bigskip

A similar result holds  with  $\partial$ instead of $t$  :

\begin{cor}
 Let $I$ be a right ideal of $A_{1}$ such that $I \cap \kk[\partial] \neq \{0\}$. Let $W$ be the p.d. subspace of $R$ such that $\theta(I) = \mathcal{D}(R,W)$. Let   $q(\partial) \in \kk[\partial]$ and $\tau : = \exp(\ad(q(\partial)))$.

Then:
\[\tau(I) = I \Longleftrightarrow \theta(q(\partial)) \in S(W) \;.\]
\end{cor}

\section{The Stafford subgroups}

Recall that if $I$ is a right ideal of $A_1$ the Stafford subgroup of $I$ is noted  $H(I)$. In the definition, the notion of isomorphism of right $A_1-$modules appears. Now, as the ring $A_1$ is hereditary, if $I,J$ are isomorphic right ideals of $A_1$, then there exists $x \in Q_1$ such that $xI =J$. So, if $I$ is  a right ideal of $A_{1}$, then we have :

i) $\qq \sigma \in \Aut_{\kk}(A_1),\; H(\sigma(I)) = \sigma H(I)\sigma\inv \;;$

ii) $\qq 0\neq z \in Q_1,\; H(zI) = H(I)$.

\bigskip

We will simply note $H(V) := H(\mathcal{D}(R,V))$ for any p.d. subspace $V$ of $R$.

Following the above remark, a Stafford subgroup of $\Aut_\kk(A_1)$ is of the form $H(V) $ for some p.d.i. subspace of $R$.

\begin{pro}\label{pro:8}

Let $V$ and $W$ be two p.d.i. subspaces of $R$.

If $H(V )\subseteq H(W)$ then  :
\[
i)\;\; S(V )\subseteq  S(W) \mbox{ and } ii)\;\; \mathcal{C}(R, V ) \subseteq \mathcal{C}(R,W)\;.\]
\end{pro}

\begin{dem}

i) : Let $p \in S(V)$. By the theorem \ref{thm:exp}, $\sigma := \exp(\ad (p))\in H(V)$, thus $\sigma\in H(W)$. So there exists $0\neq a \in Q_1$ such that $\sigma({\cal D}(R,W) = a{\cal D}(R,W)$. So $a{\cal D}(R,W) \sub A_1$. In particular, $a \in\kk(t)[\partial]$. We have also :
\[\deg_\partial a \le \deg_\partial d\]
for all $d \in {\cal D}(R,W)$. Therefore, $\deg_\partial a =0$ and $a \in \kk(t)$. But :
\[\sigma({\cal D}(R,W)) \star 1 =a {\cal D}(R,W)\star 1\]
\[= a W \]
thus $aW \sub R$ and  $a \in R$ because $RW=R$. Since $\sigma(t) =t$, we have :
\[a\inv {\cal D}(R,W) =  \sigma\inv ({\cal D}(R,W))\]
so $a\inv \in \kk[t]$ too. Therefore, $a \in\kk^*$ and 
\[\sigma({\cal D}(R,W) ) = {\cal D}(R,W)\]
\[\Longrightarrow \, p \in S(W)\]
by the theorem \ref{thm:exp}, again.

ii) : If $a \in {\cal C}(R,V)$, then :
\[a R \sub V\]
\[\Longrightarrow aRV \sub V\]
\[\Longrightarrow aR \sub S(V)\]
\[\Longrightarrow aR \sub S(W)\]
(by i))
\[\Longrightarrow aRW \sub W\]

\[\Longrightarrow aR \sub W\]
\ie $a \in {\cal C}(R,W)$.
\end{dem}

\begin{pro}\label{pro:9}
Let $I$ and $J$ be two right ideals of $A_{1}$ such that:
$\theta(I) = \mathcal{D}(R, V )$ and  $\theta(J) = \mathcal{D}(R,W)$ with $V$ and $W$  two p.d.i. subspaces of $R$. Then :
\[H(I)\subseteq H(J)\Longrightarrow  I\cap k[\partial]\subseteq J \cap\ k[\partial].\;.\]
\end{pro}
\begin{dem}
We apply the proposition \ref{pro:8} to  $\theta(I)$ and $\theta(J)$.
\end{dem}

Now, we deduce the following for the caracteristic elements of p.d. subspaces of $R$ :

\begin{cor}\label{cor:dix}

Let $V$ and $W$ be two p.d.i. subspaces of $R$ such that
$H(V)\subseteq H(W)$.

If $e^*_V \in{\cal D}(R,V)^*$ and $f_V\in{\cal D}(R,V)$ are the characteristic elements of $V$ and $e_W^*\in{\cal D}(R,W)^*$ , $f_W\in{\cal D}(R,W)$ are those of $W$, then :
\[e^*_Vf_V \in e_W^*f_W\kk[\partial]\;.\]
\end{cor}

\begin{dem}

 We have $H(\mathcal{D}(R, V )) = H(e^*_V\mathcal{D}(R, V ))$ and $H(\mathcal{D}(R,W)) = H(e_W^*\mathcal{D}(R,W))$, so we have the inclusion:
$$H(e^*_V\mathcal{D}(R, V )) \subseteq H(e_W^*\mathcal{D}(R,W))\;.$$
Now, we can show that $e^*_V$ is the unique element in $Q_1$ (modulo $\kk^*$) such that :
\[\theta (e^*_V{\cal D}(R,V)) = {\cal D}(R,V')\]
for some p.d.i. subspace $V'$ of $R$.  

By the proposition \ref{pro:9} above, we have :
\[e^*_V\mathcal{D}(R, V ) \cap \kk[\partial] \subseteq  e_W^*\mathcal{D}(R,W) \cap \kk[\partial]\;.\]

Since we have $e^*_V\mathcal{D}(R, V ) \cap \kk[\partial] = e^*_V f_V \kk[\partial]$ and $e_W^*\mathcal{D}(R,W) \cap \kk[\partial] = e_W^*f_W\kk[\partial]$, we obtain :
\[e^*_Vf_V \in e_W^*f_W\kk[\partial]\;.\]
\end{dem}

\bigskip

Now we are ready to prove the main proposition.

We will say that a p.d. subspace of $R$ is {\it monomial} if it can be generated by monomials.

\begin{pro}\label{pro:princ}
Let $V\subset R$ be a proper p.d.i. subspace of $R$ such that $H(\kk[X_n])\subseteq H(V)$. Then  :

i) $V$ is monomial ;

ii) $\mathcal{C}(R, V ) = t^{n}\kk[t]$ ;

iii) $V = \kk[X_{n}]$.

\end{pro}
\begin{dem}

 If $n = 1$, clearly $V$ would be equal to $R$, contrary to our hypothesis. So let $n \geq 2$. Let $e^*\in{\cal D}(R,V)^*$ and $f \in {\cal D}(R,V)$ be the characteristic elements of $V$.

{\bf i)} We have $\kk[X_{n}] = \kk + t^{n}\kk[t]$, $\mathcal{C}(R,\kk[X_{n}]) = t^{n}\kk [t]$ and $t^{n}\kk [t] \subseteq \mathcal{C}(R, V )$. In particular, $t^n \in {\cal D}(R,V)$ and as $e^* \in {\cal D}(R,V)^*$, we have :
\[e^*,\,f \in \kk[t,t\inv,\partial] \;.\]
 So, we can use the standard form to describe $e^*$ and $f$ :
\[ e^*= t^{p}a_{p}(t\partial) + ...+ t^{q}a_{q}(t\partial)\]
\[
f= t^{r}b_{r}(t\partial) + ... + t^{s}b_{s}(t\partial)\]
 
for some integers $p\le q$ and $r \leq s$ and some polynomials  $a_{i}(T),\,b_{j}(T) \in \kk[T]$. 

The characteristic elements of $\kk[X_n]$ are :
\[e_{n}^* := t^{-n}(t \partial)\,,\; f_{n}: = (t\partial-1)...(t\partial-(n-1))\]

and :
\[e_n^*f_n = \partial^n\;.\] 

According to the corollary \ref{cor:dix},  we have :
\[e_{n}^*f_{n} \in e^*f \kk[\partial]\]
and so \[e^{*}f = \partial^{l}\] for some integer $0\leq l\leq n$. That forces $p = q$ and $r =s$. Therefore  $V=f(R)$ is spanned by its monomial terms $t^{p+i}f_{p}(i)$, $i \ge 0$, and is monomial.

 In fact, $1 \le l \le n$. Otherwise :
\[e^*f=1 \Rightarrow e^*\in\kk(t)\]
and :
\[e^*V=e^*{\cal D}(R,V) \star 1 \sub R \]
\[\Rightarrow e^*R =e^*RV \sub R\]
\[\Rightarrow e^* \in \kk[t]\]
\[\Rightarrow e^*=f=1\]
\[\Rightarrow V=R\]
which is impossible.
\begin{center}
*
\end{center}
{\bf ii)} As $V$ is monomial and irreducible,  $1 \in V$ and so $\kk[X_{n}] \subseteq V$.

Suppose $t^{n-1}\in V$ and let us consider the automorphism $\sigma := \exp(\ad(t^{n-1}))$.  Clearly $t^{n-1}$ would belong to $S(V)$ since $V$ is monomial and $t^{n}k[t]\subseteq V$. Then we would have $\sigma\in H(V)$. By applying $\sigma$ to $H(\kk[X_n])$, we get a new inclusion:
\[\sigma H(\kk[X_n]) \sigma\inv \sub \sigma H(V)\sigma\inv \]
\begin{equation}\label{eq:incl}
\Longleftrightarrow  H(\sigma(\mathcal{D}(R,\kk[X_{n}]))) \subseteq H( V )\;.
\end{equation}
But for all $d \in A_1$, for all $j \ge 0$ :
\[\underbrace{[...[d,t^{n-1}],t^{n-1}],...,t^{n-1}]}_j = dt^{(n-1)j} - jt^{n-1}dt^{(n-1)(j-1)} \mod t^n A_1\;.\]
 
 So, for all $d \in A_1$, for all $r \in R$ :

\[\sigma(d) (r) = \sum_{j\ge 0}\frac{d(t^{(n-1)j}r)}{j!} - \sum_{j \ge 1} \frac{t^{n-1}d(t^{(n-1)(j-1)}r)}{(j-1)!} \mod t^n\kk[t]\]
(those are  in fact  finite sums because $d(t^{(n-1)j}r) \in t^n\kk[t]$
for $j>>0$).

Thus for all $d \in {\cal D}(R,\kk[X_n])$, for all $r \in R$, we get :
\[ \sigma(d)(r) = (1-t^{n-1})  \sum_{j\ge 0}\frac{d(t^{(n-1)j}r)}{j!} \mod t^n\kk[t]\]
\[\in (1-t^{n-1}) \kk[X_n] +t^n\kk[t]= \kk(1-t^{n-1}) + t^n\kk[t]\;. \]

Let $U_n :=\kk(1-t^{n-1}) + t^n\kk[t]$. We have just proved :
\[\sigma({\cal D}(R,\kk[X_n])) \sub {\cal D}(R,U_n)\;.\]
In the same way, we can prove :
\[\sigma\inv({\cal D}(R,U_n))\sub{\cal D}(R,\kk[X_n])\;.\]

Therefore, we have exactly :
\[\sigma({\cal D}(R,\kk[X_n])) = {\cal D}(R,U_n) \;.\]

Now, the characteristic elements of $U_n$ are :
\[e^*_{U_n}:= \partial ^{n-2}t^{-n}(t\partial) + (-1)^{n}(n - 1)!t^{1-n} \in \mathcal{D}(U_{n},R)\]
\[f_{U_n} = (t\partial - 1)... (t\partial - (n - 1)) + (-1)^{n}(n - 1)!t^{n-1} \in \mathcal{D}(R,U_{n}) \;.\]
We have :
\[e^*_{U_n}f_{U_n} = \left(\partial^{n-1}+(-1)^n(n-1)!\right)^2\]
hence : \[e^*_{U_n}f_{U_n} \notin \partial^l\kk[\partial]= e^*f\kk[\partial]\;.\]

By the corollary \ref{cor:dix}, we deduce :
\[H(\sigma(\mathcal{D}(R,\kk[X_{n}]))) = H(U_n) \not\subseteq H( V )\]
contrary to (\ref{eq:incl}).

So $t^{n-1} \not\in V$ and we have exactly :
\[\mathcal{C}(R, V ) = t^{n}\kk[t] \;.\]
\begin{center}
*
\end{center}
{\bf iii)} Now, $\mathcal{C}(R, V ) = t^{n}k[t]$, and $t^{n-1} \not\in V$ . For $n = 2$, we have already $V = \kk[X_{2}]$.
 So we will suppose $n\geq 3$. 

Let $1\le n_{1} < n_{2} <...< n_{s} < n - 1$ be the integers such that $t^{n_{i}}\not \in V$ . We will show that $s = n - 2$ and thus $V = \kk[X_n]$.

We use  again the automorphism $ \sigma = \exp(\ad(t^{n-1}))$. We find :

\[ \sigma(\mathcal{D}(R, V )) = \mathcal{D}(R, V_{\sigma})\]
 where $V_{\sigma} := \kk(1 - t^{n-1}) + V \cap t \kk[t]$. We set :
\[h(T) := (T-n_1)...(T-n_s) \] 
and \[\lambda := (n-1)! \frac{h(0)}{h(n-1)} \;.\]
We check that the element
\[g_{\sigma} := h(t\partial)(t\partial - (n - 1))\partial ^{n-2} + \lambda th(t\partial + 1) \]
is an element of $\mathcal{D}(R, V_{\sigma})$. We see that \[\deg_t(g_{\sigma}) = s+1\]
\[ =\dim_{\kk} {R}/{V_{\sigma}}\]
(for example because of the short exact sequence :
\[ 0 \to V_\sigma/(t\kk[t] \cap V) \to R/(V \cap t\kk[t]) \to R/V_\sigma \to 0\;).\] So $g_{\sigma}$ has minimum $t$-degree. The element  $g_{\sigma}$ can be expanded as :
\[g_{\sigma} = t^{s+1}(\partial^{n-1} + \lambda)\partial^{s} + t^{s}b_{s}(\partial) +...+ tb_{1}(\partial) - (n - 1)h(0)\partial^{n-2}\]
for some polynomials $b_i(T) \in \kk[T]$.

Since $\lambda \neq 0$, the highest common factor of 
\[(\partial ^{n-1}+\lambda)\partial ^{s}, b_{s}(\partial),... , b_{1}(\partial), \partial^{n-2}\]
must be some $\partial^{ r}$ where $0  \leq r\leq s$. So if $e^*_\sigma,f_\sigma$ are the characteristic elements of $V_\sigma$, we have $g_{\sigma} = f_{\sigma}\partial^{r}$. Thus  $f_{\sigma}$
must be equal to :
\[t^{s+1}(\partial ^{n-1} + \lambda)\partial^{s-r} + t^{s}a_{s}(\partial) +... + ta_{1}(\partial) - (n - 1)h(0)\partial ^{n-2-r}\]
where $a_i(\partial) := b_i(\partial)\partial^{-r} \in \kk[\partial]$.

But :
\[H(\kk[X_n]) \sub H(V)\]
\[\Longrightarrow \sigma H(\kk[X_n])\sigma\inv \sub \sigma H(V)\sigma\inv\]
\[\Longrightarrow H(U_n) \sub H(V_\sigma)\]
\[\Longrightarrow e^*_{U_n}f_{U_n} \in e^*_\sigma f_\sigma\kk[\partial]\]
by the corollary \ref{cor:dix}. Now, by the proposition \ref{pro:cinq}.iii), \[e_\sigma^*f_\sigma \in (\partial^{n-1} + \lambda)\partial^{s-r}\kk[\partial]\;.\] As a  consequence :
\[\left(\partial^{n-1} +(-1)^n(n-1)!\right)^2 \in (\partial^{n-1} + \lambda)\partial^{s-r}\kk[\partial]\]
which implies that $s=r$ and $\lambda=(-1)^n (n-1)!^2$.

Since $s=r$, we have :
\[g_\sigma =f_\sigma\partial^s\]
and :\[ g_\sigma(t^i)=0\]
for all $0 \le i \le s-1$. Thus $h(i+1)=0$ for $i=0,1,...,s-1$. But $n_1,...,n_s$ are the only roots of  $h$ so :
\[\qq 0 \le j \le s,\, n_j = j\;.\]  

From the equality $\lambda= (-1)^n (n-1)!$ we then deduce:
\[\frac{(-1)^s n_1...n_s}{(n-2)...(n-1-s)} = (- 1)^n\]
\[\equi \frac{s!(n-2-s)!}{(n-2)!} =(-1)^{n+s}\]
\[\equi { {n-2} \choose s} = (-1)^{n+s}\]
\[\equi s= 0 \mbox{ and $n$ is even  or } s= n-2 \;.\]

If $s=n-2$, then $V= \kk[X_n]$ and the proof is finished. If $s=0$ and $n$ is even, then :
\[V = \kk +\kk t +... +\kk t^{n-2}+t^n\kk[t] \]
with $n \ge 4$.

We set $\sigma:=\exp(\ad(t^{n-2}))$. Then we have :
\[H(\kk[X_n])\sub H(V)\]
\[\equi \sigma H(\kk[X_n])\sigma\inv \sub \sigma H(V) \sigma\inv \]
\[\equi H(W_n) \sub H(V_\sigma)\]
where :
\[W_n = (1-t^{n-2})\kk[X_n] +t^n\kk[t] \]
\[= \kk(1-t^{n-2}) + t^n\kk[t]\]
and \[ V_\sigma = (1-t^{n-2})V + t^n\kk[t]\]
\[= \kk + \kk(t-t^{n-1}) +\kk t^2 +...+ \kk t^{n-2} +t^n\kk [t] \;.\]

Now, let $e_{W_n}^*,f_{W_n}$ be the characteristic elements of $W_n$ and $e_\sigma^*,f_\sigma$ those of $V_\sigma$. Because of the corollary \ref{cor:dix}, we should have :

\begin{equation}\label{eq:der}
e_{W_n}^*f_{W_n} \in e_\sigma^*f_\sigma \kk[\partial] \;.
\end{equation}
But we can check that :
\[e_{W_n}^* = \left( \frac{\partial^{n-1}}{(n-1)!} + (-1)^{n-1} \partial\right) t^{1-n} + \left( \frac{\partial^{n-2}}{(n-2)!} + (-1)^{n-1} \right) t^{-n}\;,\]
\[f_{W_n} = \frac{(t\partial -1)...(t\partial -(n-1))}{(n-1)!} +(-1)^{n-1}t^{n-2}(t\partial -1) \;,\]
\[f_\sigma = t\left(\frac{\partial^{n-1}}{(n-1)!} + \partial \right) - \frac{\partial^{n-2}}{(n-2)!} -1 \;.\]

We deduce that :
\[e_{W_n}^* f_{W_n} = \left( \frac{\partial^{n-1}}{(n-1)!} + (-1)^{n-1} \partial  \right)^2\] 
\[ = \left( \frac{\partial^{n-1}}{(n-1)!}  - \partial  \right)^2\]
because $n$ is even and :
\[e_\sigma^*f_\sigma \kk[\partial] \sub  \left( \frac{\partial^{n-1}}{(n-1)!}  + \partial  \right) \kk[\partial]\]

which contradicts (\ref{eq:der}). 

Hence $V=\kk[X_n]$.
 \end{dem}

Using the description of right ideals of $A_1$ in \cite{kou}, we deduce the following :

\begin{cor}

For any non principal right ideals $I$ and $J$
, the following equivalences are satisfied:
\[H(I) \subset H(J)\]
\[\Longleftrightarrow  H(I) = H(J)\]
\[\Longleftrightarrow {}^\exists x \in  \Frac(A_{1}),\,{}^\exists \sigma \in \Aut_\kk (A_1),\;  I = x\sigma(J)\;.\]

\end{cor}

We now obtain the announced result.

\begin{pro}

Let $I$ be any right ideal of $A_{1}$. The subgroup $H(I)$ is equal to its own
normalizer subgroup in $Aut_{\kk}(A_{1})$.
\end{pro}
\bigskip
\begin{dem}
 By \cite{kou}, we can suppose $I = \mathcal{D}(R,\kk[X_n])$. Let
$\gamma \in Aut(A_{1})$ such that
\[\gamma H(\kk[X_{n}])\gamma\inv  = H(\kk[X_{n}])\;.\]

Then we have \[H(\gamma(\mathcal{D}(R,\kk[X_n])) =H(\kk[X_{n}])\;.\]

We have also : $\gamma(\mathcal{D}(R,\kk[X_n])) \simeq {\cal D}(R,V) $ for some p.d.i. subspace $V$ of $R$.

Thus, we  get $H(\kk[X_n])=H(V)$, so $V=\kk[X_n]$ by the proposition \ref{pro:princ}. Finally, we have $\gamma(\mathcal{D}(R,\kk[X_{n}])) \simeq \mathcal{D}(R,\kk[X_{n}])$ and that means $\gamma \in H(I)$.

\end{dem}
\bigskip
\bibliographystyle{plain}
\bibliography{biblio}

\end{document}